# A NOTE ON THE FRACTIONAL PERIMETER AND INTERPOLATION

AUGUSTO C. PONCE AND DANIEL SPECTOR

ABSTRACT. We present the fractional perimeter as a set-function interpolation between the Lebesgue measure and the perimeter in the sense of De Giorgi. Our motivation comes from a new fractional Boxing inequality that relates the fractional perimeter and the Hausdorff content and implies several known inequalities involving the Gagliardo seminorm of the Sobolev spaces $W^{\alpha,1}$ of order $0 < \alpha < 1$.

## 1. INTRODUCTION AND MAIN RESULTS

In the forthcoming work [10], the authors established a family of estimates that in some sense interpolate the classical Boxing inequality of W. Gustin [6] and a trivial endpoint for the Lebesgue measure:

**Theorem 1.** *There exists a constant $C > 0$ depending only on the dimension $d \geq 1$ such that*

$$\mathcal{H}_\infty^{d-\alpha}(U) \leq C\alpha(1-\alpha) P_\alpha(U),$$

*for every bounded open subset $U \subset \mathbb{R}^d$, uniformly with respect to $\alpha \in (0,1)$. Here, for a measurable set $A \subset \mathbb{R}^d$,*

$$(1.1) \qquad \mathcal{H}_\infty^{d-\alpha}(A) := \inf\left\{ \sum_{i=0}^\infty \omega_{d-\alpha} r_i^{d-\alpha} : A \subset \bigcup_{i=0}^\infty B(x_i, r_i) \right\}$$

*is the Hausdorff content of dimension $(d-\alpha)$ and*

$$(1.2) \qquad P_\alpha(A) := 2 \int_A \int_{\mathbb{R}^d \setminus A} \frac{\mathrm{d}y \, \mathrm{d}z}{|y-z|^{\alpha+d}}$$

*is the fractional perimeter of $A$.*

When combined with a straightforward counterpart of the coarea formula for the fractional Sobolev space $W^{\alpha,1}(\mathbb{R}^d)$, one obtains the following trace inequality for every continuous function $u \in W^{\alpha,1}(\mathbb{R}^d)$:

$$(1.3) \qquad \int_{\mathbb{R}^d} |u| \, \mathrm{d}\mu \leq C\alpha(1-\alpha) \int_{\mathbb{R}^d} \int_{\mathbb{R}^d} \frac{|u(x)-u(y)|}{|x-y|^{\alpha+d}} \, \mathrm{d}y \, \mathrm{d}x,$$

where $\mu$ is any nonnegative Borel measure in $\mathbb{R}^d$ such that $\mu(B(x,r)) \leq \omega_{d-\alpha} r^{d-\alpha}$ for all balls $B(x,r)$. This estimate is a strong form of Sobolev's inequality which implies the classical embedding of $W^{\alpha,1}$ into $L^{\frac{d}{d-\alpha}}$, its Lorentz-space improvement, and also Hardy's inequality.

The proof in [10] does not utilize interpolation, and in fact the authors do not know of any framework that allows one to interpolate the Hausdorff content. The fractional perimeter, however, is often thought of as an object which is intermediate between the perimeter (in the sense of De Giorgi) and the Lebesgue measure. For example, one has the fractional Gagliardo-Nirenberg inequality inherited from that for $BV$ functions (see e.g. Proposition 4.2 in [2] or Proposition 15.6 in [9]):

$$\text{(1.4)} \qquad \alpha(1-\alpha) P_\alpha(U) \leq C \, |U|^{1-\alpha} \operatorname{Per}(U)^\alpha,$$

for all $\alpha \in (0,1)$, as well as the asymptotics

$$\text{(1.5)} \qquad \lim_{\alpha \to 0} \alpha P_\alpha(U) = C'|U| \quad \text{and} \quad \lim_{\alpha \to 1} (1-\alpha) P_\alpha(U) = C'' \operatorname{Per}(U),$$

that allows recovery of the endpoints; see [3, 5, 7].

The purpose of this note is to give a sense in which the fractional perimeter is intermediate between the perimeter and the Lebesgue measure. To this end, let us recall some results that follow from a routine application of the real interpolation theory. First, one has

$$P_\alpha(A) = [\chi_A]_{W^{\alpha,1}(\mathbb{R}^d)} := \int_{\mathbb{R}^d} \int_{\mathbb{R}^d} \frac{|\chi_A(y) - \chi_A(z)|}{|y-z|^{\alpha+d}} \, dy \, dz,$$

that is, the fractional perimeter of a measurable set $A \subset \mathbb{R}^d$ is the Gagliardo semi-norm on $W^{\alpha,1}(\mathbb{R}^d)$ applied to the function $\chi_A$. The fractional Sobolev space $W^{\alpha,1}(\mathbb{R}^d)$ is itself a Besov space which arises in the real interpolation of $L^1(\mathbb{R}^d)$ and $\dot{W}^{1,1}(\mathbb{R}^d)$ with parameter $\alpha$ (see e.g. Corollary 4.13 and Eq. (4.42) in [1]):

$$\dot{W}^{\alpha,1}(\mathbb{R}^d) = (L^1, \dot{W}^{1,1})_\alpha,$$

though in the real interpolation method one encounters a variety of equivalent semi-norms. Indeed, a convenient method for computing the semi-norm of a function $f \in (L^1, \dot{W}^{1,1})_\alpha$ is to introduce the $K$-functional

$$K(t, f, L^1, \dot{W}^{1,1}) := \inf_{f = f_1 + f_2} \|f_1\|_{L^1(\mathbb{R}^d)} + t \|\nabla f_2\|_{L^1(\mathbb{R}^d)},$$

where $t > 0$ and the infimum is taken over all decompositions $f = f_1 + f_2$ such that $f_1 \in L^1(\mathbb{R}^d)$ and $f_2 \in \dot{W}^{1,1}(\mathbb{R}^d)$. One then obtains a semi-norm on the interpolation space $(L^1, \dot{W}^{1,1})_\alpha$ via the formula

$$[f]_{\tilde{W}^{\alpha,1}(\mathbb{R}^d)} := \alpha(1-\alpha) \int_0^\infty t^{-\alpha} K(t, f, L^1, \dot{W}^{1,1}) \, \frac{dt}{t},$$

for which one can show

$$[f]_{\tilde{W}^{\alpha,1}(\mathbb{R}^d)} \sim \alpha(1-\alpha) [f]_{W^{\alpha,1}(\mathbb{R}^d)};$$

see [8]. Here, the symbol $\sim$ indicates that both quantities are comparable, uniformly with respect to $\alpha \in (0,1)$. While this is a standard approach to $\dot{W}^{\alpha,1}(\mathbb{R}^d)$, it has the defect of not being applicable to the characteristic functions of sets, even sets of finite perimeter, to deduce estimate (1.4) and some type of analogue to the limits (1.5). This is not a serious setback, as it is not difficult to show the equivalence

$$K(t, f, L^1, \dot{W}^{1,1}) = K(t, f, L^1, \dot{BV}),$$

which is connected to the assertion that the Gagliardo closure of $\dot{W}^{1,1}$ is $\dot{BV}$. This can be deduced from the literature, for example, through an application of Lemma 2 of the paper of Cwikel [4, pp. 216-217] and as a result one obtains

(1.6) $$\dot{W}^{\alpha,1}(\mathbb{R}^d) = (L^1, \dot{BV})_\alpha.$$

This yields Eq. (1.4) directly from the classical theory, while an analogue to Eq. (1.5),

$$\lim_{\alpha \to 0} [\chi_A]_{\tilde{W}^{\alpha,1}(\mathbb{R}^d)} = |A| \quad \text{and} \quad \lim_{\alpha \to 1} [\chi_A]_{\tilde{W}^{\alpha,1}(\mathbb{R}^d)} = \operatorname{Per}(A),$$

follows (cf. [8]) from the fact that the pair $(L^1(\mathbb{R}^d), \dot{BV}(\mathbb{R}^d))$ is normal, i.e.

$$\lim_{t \to \infty} K(t, f, L^1, \dot{BV}) = \|f\|_{L^1(\mathbb{R}^d)} \quad \text{and} \quad \lim_{t \to 0} \frac{K(t, f, L^1, \dot{BV})}{t} = |Df|(\mathbb{R}^d).$$

These results motivate us to define a fractional perimeter intrinsic to interpolation, which is

$$\tilde{P}_\alpha(A) := [\chi_A]_{\tilde{W}^{\alpha,1}(\mathbb{R}^d)} = \alpha(1-\alpha) \int_0^\infty t^{-\alpha} K(t, \chi_A, L^1, \dot{BV}) \frac{dt}{t},$$

for which standard interpolation arguments yield the following refinements to Eq. (1.4) and (1.5) (whose proofs we supply for the convenience of the reader).

**Theorem 2.** *If $A \subset \mathbb{R}^d$ is a set of finite perimeter, then one has*

$$\tilde{P}_\alpha(A) \leq |A|^{1-\alpha} \operatorname{Per}(A)^\alpha,$$

*and*

$$\lim_{\alpha \to 0} \tilde{P}_\alpha(A) = |A| \quad \text{and} \quad \lim_{\alpha \to 1} \tilde{P}_\alpha(A) = \operatorname{Per}(A).$$

Yet this approach still views sets of fractional finite perimeter as intermediate between $L^1$ and $\dot{BV}$, while one might wish to use directly the spaces of sets of finite Lebesgue measure and those of finite perimeter. We here pursue this approach, though as these are not linear spaces we prefer to do so via the penalty functional

$$K(t, A) = \inf_{U \subset \mathbb{R}^d} |A \triangle U| + t \operatorname{Per}(U),$$

where the infimum is taken over all open sets $U \subset \mathbb{R}^d$ of finite perimeter and $A \triangle U := (A \setminus U) \cup (U \setminus A)$ is the symmetric difference between the sets $A$ and $U$. Notice that

$$K(t, \chi_A, L^1, \dot{BV}) \le K(t, A),$$

and so we might hope to use $K(t, A)$ directly in the computation of the fractional perimeter. Indeed, we have

**Theorem 3.** *For any Lebesgue measurable set $A \subset \mathbb{R}^d$ and $t > 0$, one has*

$$K(t, \chi_A, L^1, \dot{BV}) = K(t, A)$$

*and so in particular*

$$\tilde{P}_\alpha(A) = \alpha(1-\alpha) \int_0^\infty t^{-\alpha} K(t, A) \frac{dt}{t}.$$

Moreover, the argument in interpolation that demonstrates the equivalence between the semi-norms $[f]_{W^{\alpha,1}}$ and $[f]_{\tilde{W}^{\alpha,1}}$ easily gives the following geometric interpretation of the penalization functional $K$ (see Theorem 4.12 and (4.42) in [1]):

**Theorem 4.** *For any Lebesgue measurable set $A \subset \mathbb{R}^d$ and $t > 0$, one has*

$$K(t, A) \sim \sup_{h \in B(0,t)} |(A + h) \triangle A|.$$

When one recalls the elementary comparison

$$\sup_{h \in B(0,t)} |(A + h) \triangle A| \sim \fint_{B(0,t)} |(A+h) \triangle A| \, dh,$$

we find

$$\tilde{P}_\alpha(A) \sim \alpha(1-\alpha) \int_0^\infty t^{-\alpha} \left( \fint_{B(0,t)} |(A+h) \triangle A| \, dh \right) \frac{dt}{t}$$

and then Fubini's theorem yields $\tilde{P}_\alpha(A) \sim \alpha(1-\alpha) P_\alpha(A)$. This completes our work on describing the space of sets of fractional finite perimeter as being between the space of sets of finite Lebesgue measure and the space of sets of finite perimeter.

## 2. Proofs of the Main Results

*Proof of Theorem 2.* Given a set with finite perimeter $A \subset \mathbb{R}^d$, for every $t > 0$ we have

$$K(t, A) \le t \operatorname{Per}(A) \quad \text{and} \quad K(t, A) \le |A|.$$

Thus, for $r > 0$ to be explicitly chosen later on,

$$\int_0^\infty t^{-\alpha} K(t,A) \, \frac{\mathrm{d}t}{t} = \int_0^r t^{-\alpha} K(t,A) \, \frac{\mathrm{d}t}{t} + \int_r^\infty t^{-\alpha} K(t,A) \, \frac{\mathrm{d}t}{t}$$

$$\leq \operatorname{Per}(A) \int_0^r t^{-\alpha} \, \mathrm{d}t + |A| \int_r^\infty t^{-1-\alpha} \, \mathrm{d}t$$

$$= \operatorname{Per}(A) \frac{r^{1-\alpha}}{1-\alpha} + |A| \frac{r^{-\alpha}}{\alpha}.$$

Minimizing the right-hand side with respect to $r$ we get

$$\int_0^\infty t^{-\alpha} K(t,A) \, \frac{\mathrm{d}t}{t} \leq \frac{1}{\alpha(1-\alpha)} |A|^{1-\alpha} \operatorname{Per}(A)^\alpha.$$

To conclude, it suffices to observe that

$$\lim_{t \to \infty} K(t,A) = |A| \quad \text{and} \quad \lim_{t \to 0} \frac{K(t,A)}{t} = \operatorname{Per}(A),$$

which immediately implies the limits for $\tilde{P}_\alpha(A)$. □

*Proof of Theorem 3.* Observe that $\|\chi_A - \chi_U\|_{L^1(\mathbb{R}^d)} = |A \triangle U|$ and $|D\chi_U|(\mathbb{R}^d) = \operatorname{Per}(U)$. The decomposition $\chi_A = (\chi_A - \chi_U) + \chi_U$ thus yields the inequality

$$K(t, \chi_A, L^1, \dot{BV}) \leq K(t, A),$$

and so we must prove the reverse inequality.

For every $f \in \dot{BV}(\mathbb{R}^d)$ and $t > 0$, we begin by showing that

$$(2.1) \qquad K(t, f, L^1, \dot{BV}) = \inf_{\substack{f = f_1 + f_2 \\ f_2 \in C_c^\infty(\mathbb{R}^d)}} \|f_1\|_{L^1(\mathbb{R}^d)} + t \|Df_2\|_{L^1(\mathbb{R}^d)}.$$

The inequality $\leq$ is immediate since the infimum in the left-hand side is taken over a larger class of decompositions of $f$. For the reverse inequality, we let $\eta > 0$ and $f = f_1 + f_2$ be such that $f_2 \in \dot{BV}(\mathbb{R}^d)$ and

$$\|f_1\|_{L^1(\mathbb{R}^d)} + t |Df_2|(\mathbb{R}^d) \leq K(t, f, L^1, \dot{BV}) + \eta.$$

Take $g_2 \in C_c^\infty(\mathbb{R}^d)$ such that

$$\|f_2 - g_2\|_{L^1(\mathbb{R}^d)} \leq \eta \quad \text{and} \quad |Dg_2|(\mathbb{R}^d) \leq |Df_2|(\mathbb{R}^d) + \frac{\eta}{t}.$$

The decomposition $f = (f_1 + f_2 - g_2) + g_2$ satisfies

$$\|f_1 + f_2 - g_2\|_{L^1(\mathbb{R}^d)} + t|Dg_2|(\mathbb{R}^d)$$
$$\leq \|f_1\|_{L^1(\mathbb{R}^d)} + \|f_2 - g_2\|_{L^1(\mathbb{R}^d)} + t|Dg_2|(\mathbb{R}^d)$$
$$\leq K(t, f, L^1, \dot{BV}) + 3\eta.$$

Since $\eta > 0$ is arbitrary, this concludes the proof of the identity (2.1).

By Eq. (2.1), we may thus restrict our attention to decompositions $\chi_A = (\chi_A - g) + g$ with $g \in C_c^\infty(\mathbb{R}^d)$. We first observe that the choice $U := \{|g| > s\}$ with $s \in (0,1)$ satisfies

$$(2.2) \qquad A \triangle U \subset \big(\{|\chi_A - g| \geq s\} \setminus A\big) \cup \big(\{|\chi_A - g| \geq 1 - s\} \cap A\big).$$

We now explain how to make a suitable choice of $s$. To this end, by Cavalieri's principle and the classical coarea formula we have

$$\int_0^\infty |\{|\chi_A - g| \geq s\}|\,\mathrm{d}s + t\int_0^\infty \mathrm{Per}(\{|g| > s\})\,\mathrm{d}s \\ = \|\chi_A - g\|_{L^1(\mathbb{R}^d)} + t\|Dg\|_{L^1(\mathbb{R}^d)}.$$

By comparison of integrals and an affine change of variables we have

$$\int_0^\infty |\{|\chi_A - g| \geq s\}|\,\mathrm{d}s \\ \geq \int_0^1 |\{|\chi_A - g| \geq s\}|\,\mathrm{d}s \\ = \int_0^1 \Big(|\{|\chi_A - g| \geq s\} \setminus A| + |\{|\chi_A - g| \geq 1 - s\} \cap A|\Big)\,\mathrm{d}s.$$

Thus,

$$\int_0^1 \Big(|\{|\chi_A - g| \geq s\} \setminus A| + |\{|\chi_A - g| \geq 1 - s\} \cap A| + t\,\mathrm{Per}\,(\{|g| > s\})\Big)\,\mathrm{d}s \\ \leq \|\chi_A - g\|_{L^1(\mathbb{R}^d)} + t\|Dg\|_{L^1(\mathbb{R}^d)}.$$

Take $s \in (0,1)$, depending on $t$, such that the open set $\{|g| > s\}$ is smooth and

$$|\{|\chi_A - g| \geq s\} \setminus A| + |\{|\chi_A - g| \geq 1 - s\} \cap A| + t\,\mathrm{Per}\,(\{|g| > s\}) \\ \leq \|\chi_A - g\|_{L^1(\mathbb{R}^d)} + t\|Dg\|_{L^1(\mathbb{R}^d)}.$$

In view of (2.2) we get

$$K(t, A) \leq |A \triangle U| + t\,\mathrm{Per}\,(U) \leq \|\chi_A - g\|_{L^1(\mathbb{R}^d)} + t\|Dg\|_{L^1(\mathbb{R}^d)}.$$

Taking the infimum of the right-hand side with respect to $g$ we deduce the reverse inequality. $\square$

*Proof of Theorem 4.* Let $h \in B(0,t)$. For every subset $U \subset \mathbb{R}^d$ of finite perimeter, by the triangle inequality in $L^1$ and the translation invariance of the Lebesgue measure we have

$$|(A + h) \triangle A| = \|\chi_{A+h} - \chi_A\|_{L^1(\mathbb{R}^d)} \\ \leq \|\chi_{A+h} - \chi_{U+h}\|_{L^1(\mathbb{R}^d)} + \|\chi_A - \chi_U\|_{L^1(\mathbb{R}^d)} + \|\chi_{U+h} - \chi_U\|_{L^1(\mathbb{R}^d)} \\ = 2|A \triangle U| + |(U + h) \triangle U|.$$

Since $U$ has finite perimeter, $|(U+h) \triangle U| \leq C_1 |h|\,\mathrm{Per}\,(U)$, and since $|h| \leq t$ this implies

$$|(A+h) \triangle A| \leq 2|A \triangle U| + C_1|h|\,\mathrm{Per}\,(U) \leq C_2\,K(t, A).$$

To get the reverse comparison, take a smooth mollifier $\rho$ supported in the ball $B(0,t)$ and write

$$\chi_A = (\rho * \chi_A - \chi_A) + \rho * \chi_A.$$

Observe that

$$\|\rho * \chi_A - \chi_A\|_{L^1(\mathbb{R}^d)} \leq \sup_{h \in B(0,t)} |(A+h) \triangle A|,$$

while $\int_{\mathbb{R}^d} D\rho = 0$ and the fact that we can choose $\rho$ such that $\|D\rho\|_{L^\infty(\mathbb{R}^d)} \leq C_3/t$ implies

$$\|D(\rho * \chi_A)\|_{L^1(\mathbb{R}^d)} \leq \int_{\mathbb{R}^d} \int_{\mathbb{R}^d} |D\rho(h)|[\chi_{A+h}(x) - \chi_A(x)] \, \mathrm{d}h \, \mathrm{d}x$$
$$\leq \frac{C_3}{t} \sup_{h \in B(0,t)} |(A+h) \triangle A|.$$

It thus remains to argue as in the proof of Theorem 3 and take $U = \{\rho * \chi_A > s\}$ for some suitable $s \in (0,1)$. □

## ACKNOWLEDGEMENTS

The authors would like to warmly thank Mario Milman for helpful discussions, comments, and references. D.S. is supported by the Taiwan Ministry of Science and Technology under research grant number 105-2115-M-009-004-MY2.

AUGUSTO C. PONCE
UNIVERSITÉ CATHOLIQUE DE LOUVAIN
INSTITUT DE RECHERCHE EN MATHÉMATIQUE ET PHYSIQUE
CHEMIN DU CYCLOTRON 2, L7.01.02
1348 LOUVAIN-LA-NEUVE, BELGIUM
*E-mail address*: `Augusto.Ponce@uclouvain.be`

DANIEL SPECTOR
NATIONAL CHIAO TUNG UNIVERSITY
DEPARTMENT OF APPLIED MATHEMATICS
HSINCHU, TAIWAN

NATIONAL CENTER FOR THEORETICAL SCIENCES
NATIONAL TAIWAN UNIVERSITY
NO. 1 SEC. 4 ROOSEVELT RD.
TAIPEI, 106, TAIWAN
*E-mail address*: `dspector@math.nctu.edu.tw`